\documentclass[onefignum,onetabnum]{nosiam}

\usepackage{lipsum}
\usepackage{amsfonts}
\usepackage{graphicx}
\usepackage{epstopdf}
\usepackage{algorithmic}
\Crefname{ALC@unique}{Line}{Lines}
\ifpdf
  \DeclareGraphicsExtensions{.eps,.pdf,.png,.jpg}
\else
  \DeclareGraphicsExtensions{.eps}
\fi

\usepackage[utf8]{inputenc}
\usepackage{geometry}
\usepackage{amsmath, amssymb} 
\hypersetup{
    colorlinks=true,
    linkcolor=black,
    filecolor=black,      
    urlcolor=black,
    citecolor=black
} 

\newtheorem{proposition}{Proposition}

\newtheorem{definition}{Definition}
  \theoremstyle{plain}
  \theoremheaderfont{\normalfont\itshape}
  \theorembodyfont{\normalfont}
\newtheorem{remark}{Remark}

\usepackage{enumitem}
\setlist[enumerate]{leftmargin=.5in}
\setlist[itemize]{leftmargin=.5in}

\newsiamremark{hypothesis}{Hypothesis}
\crefname{hypothesis}{Hypothesis}{Hypotheses}
\newsiamthm{claim}{Claim}

\headers{Robust geometric modeling of 3-periodic tensegrity frameworks}{M. Himmelmann, and M. E. Evans}

\title{Robust geometric modeling of 3-periodic tensegrity frameworks using Riemannian optimization\thanks{Funded by the Deutsche Forschungsgemeinschaft (DFG - German Research Foundation) - Project-ID 195170736 - TRR109}}

\author{Matthias Himmelmann\thanks{University of Potsdam, Institute for Mathematics, Karl-Liebknecht-Str. 24-25, 14476 Potsdam, Germany 
  (\email{himmelmann1@uni-potsdam.de}).}
\and Myfanwy E. Evans\thanks{University of Potsdam, Institute for Mathematics, Karl-Liebknecht-Str. 24-25, 14476 Potsdam, Germany 
  (\email{evans@uni-potsdam.de}).}}

\usepackage{amsopn}

\makeatletter
\newcommand*{\addFileDependency}[1]{
  \typeout{(#1)}
  \@addtofilelist{#1}
  \IfFileExists{#1}{}{\typeout{No file #1.}}
}
\makeatother

\begin{document}

\maketitle

\begin{abstract}
  Framework materials and their deformations provide a compelling relation between materials science and algebraic geometry. Physical distance constraints within the material transform into polynomial constraints, making algebraic geometry and associated numerical strategies feasible for finding equilibrium configurations and deformation pathways. In this paper, we build the necessary geometric formulations and numerical strategies to explore the mechanics of two examples of 3-periodic tensegrity frameworks through non-linear optimization, eventually showing that the structures are auxetic by multiple definitions. 
\end{abstract}



\section{Introduction}

In a geometric context, framework materials are structures composed of many rigid bars (or edges), pinned together at vertices where the bars can freely rotate relative to each other. These structures can be rigid or deformable, depending on the combinatorics of how the bars are joined together. For example, take three bars of equal length connected to form a triangle, or four bars connected into a square; the triangle is rigid, whereas the square can be deformed into a family of parallelograms. With increasing complexity, in particular periodicity, frameworks can be used as a simplified model of various materials, such as zeolites and metal-organic frameworks~\cite{MOK_rods2}. Due to their high porosity, mechanics and other properties, which stem from the underlying geometric structure, framework materials are engineered for various applications. The precise understanding of geometric structure in framework materials is thus a valuable contribution to enhanced material design.

Possible deformation pathways of a framework structure can be described algebraically. Each rigid rod can be represented by a polynomial constraint on its end points; they are always separated by a distance equal to the length of the rod. The solution set of the polynomial system arising when accumulating all such constraints provides possible configurations of the framework, be it isolated real solutions for a rigid structure, or one- or higher-dimensional configuration spaces for flexible structures. There is thus a clear connection between framework behavior and solutions to systems of polynomials. 

Auxetic frameworks, or auxetic materials more generally, can be characterized as those structures that exhibit a perpendicular expansion upon stretching the material in a chosen direction; a somewhat counterintuitive property. This is quantified by a negative Poisson's ratio~\cite{Lakes1987}. Such materials can occur naturally, but are also a prime target in the geometric design of metamaterials with prescribed microstructures and targeted functionality, such as impact protection and filtration. In a mathematical context, Borcea and Streinu have introduced a geometric description of auxetic deformations \cite{geometricauxetics, Borcea2019AuxeticRI}, which we use extensively in parts of this article. 

\begin{figure}[h!]
    \centering
    \includegraphics[width=0.33\linewidth]{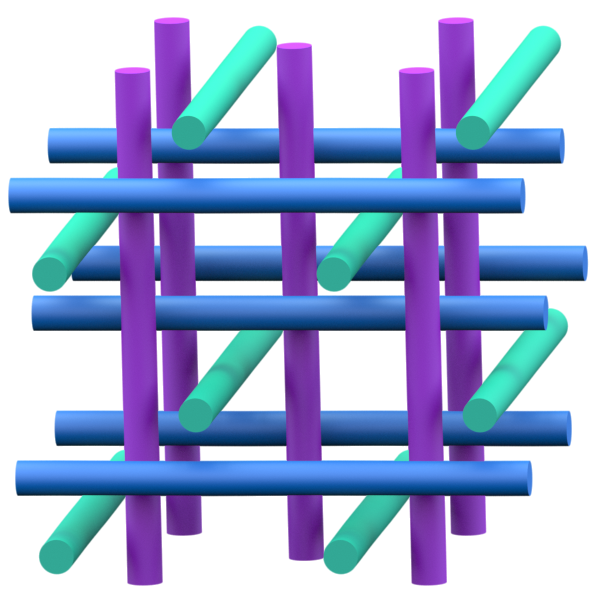}
    \hspace*{1.4cm}
    \includegraphics[width=0.33\linewidth]{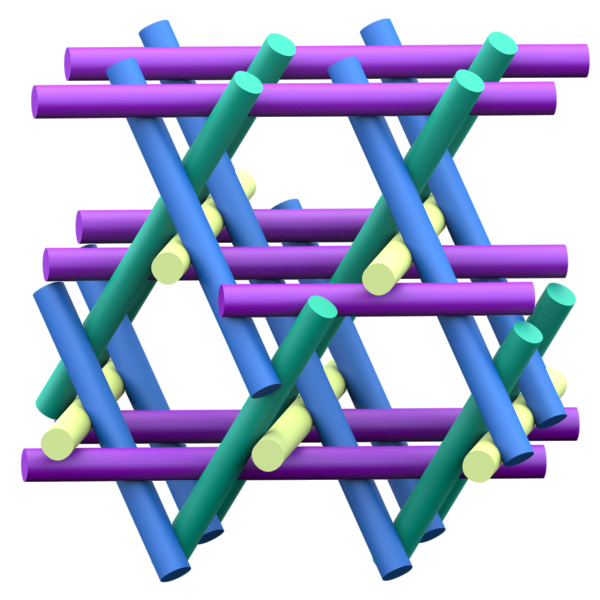}
    \caption{\small The chiral rod packings $\Pi^+$ (l.) and $\Sigma^+$ (r.), which have rods aligned along 3 or 4 axes respectively, each colored differently.}
    \label{fig:bmnsgn}
\end{figure}

Our target in this paper is the analysis of a set of slightly more complicated structures than a simple framework, namely tensegrity structures derived from filament packings. We start with a particular family of three-periodic curvilinear cylinder packings, discussed extensively in the structural chemistry literature in relation to metal-organic frameworks~\cite{MOK_rods2}. The basic geometry consists of cubic rod packings, where rods lie along the invariant axes of the cubic crystallographic symmetry groups~\cite{MOK_rods}. These packings have also been explored in a geometric setting, where their curved counterparts induced an interesting dilatant property in the material~\cite{Evans2011,Evans2015}. Here, the curved (mostly helical) cylinders have a cooperative unwinding mechanism that expands the material isotropically, reminiscent of an auxetic deformation, driven purely by the geometry of the cylinders. In a biological context, this expansion mechanism has been proposed as the strategy by which human skin cells expand and imbibe water like a sponge, driven purely by geometry~\cite{evans2014b}. The keratin intermediate filaments of the cell microstructure form one of the cubic rod packings well-known in structural chemistry~\cite{MOK_rods}. \Cref{fig:bmnsgn} shows two such rod packings, called the $\Pi^+$ packing and the $\Sigma^+$ packing, with the latter being the packing seen in the skin structure.

In a previous study, the dilatant property of the $\Pi^+$ filament packings was considered through a 3-periodic tensegrity structure~\cite{reentranttensegrity}, designed to represent the packing constraints of the original structure. Tensegrity structures are like framework materials, except some of the rigid bars are replaced by elastic elements under tension. Deformation of the tensegrity structure corresponding to the $\Pi^+$ rod packing, which is dilatant, demonstrated auxetic behavior, albeit with somewhat unstable numerics. Given this foundation, we will perform an enhanced investigation of these tensegrity structures, where the strategy for designing the tensegrity is improved for stability in \Cref{section:robustmodel}, as well as analyzing another dilatant rod packing, $\Sigma^+$, which is the skin packing described above.

In addition to these geometric considerations, a robust algorithm plays a crucial role in modeling the behavior of filament packings. Considering the framework's elastic elements as one-sided Hookean springs, its equilibration becomes a nonlinear optimization problem. As all of the proposed constraints can locally be expressed as polynomials, the constraint set can be interpreted as an algebraic variety. In \Cref{section:riemannianoptimization}, we investigate how to develop a suitable Riemannian optimization algorithm. Given the cylinder packing's contact graph as initialization, the Euclidean distance retraction \cite{euclideandistanceretraction} enables Riemannian gradient descent \cite{boumal2020intromanifolds} to yield a static energy minimum for a given initial lattice extension. Afterwards, parametrizing a path in the deformation space by one fixed lattice direction gives rise to a linear homotopy describing the system's quasistatic configurations. By varying this parameter, stress is induced on the framework and we track the resulting behavior using homotopy continuation \cite{Wampler2013bertini, HC.jl, Hom4PSArticle}. 

As described above, we will consider this enhanced analysis of the $\Pi^+$ and $\Sigma^+$ packings. We begin with the necessary background on geometric auxetics and periodic tensegrity frameworks. 

\section{Periodic Frameworks and Auxeticity}
\label{section:geometricauxetics}
We introduce here the general concepts of periodic graphs and periodic frameworks, with the goal of defining auxeticity mathematically. This will provide the geometric perspective on analyzing auxetic deformations of more complicated structures. The ideas of Definitions \ref{def:periodicgraph} and \ref{def:oneparameterdeformation} appear in \cite{Borcea2019AuxeticRI}.
\begin{definition}
\label{def:periodicgraph}
A $d$-periodic graph $(G,\Gamma)$ is a simple, connected, infinite graph $G=(V,E)$ with finite degree at every vertex together with a free abelian periodicity group $\Gamma\subset \text{Aut}(G)$ of rank $d$, with finitely many vertex orbits $V/\Gamma$ and edge orbits $E/\Gamma$.

We define a $d$-periodic framework $\mathcal{F}=(G,\Gamma, p,\pi)$ by a placement $p:V\rightarrow \mathbb{R}^d$ of the $d$-periodic graph $G$'s vertices and an injective group representation $\pi:\Gamma \rightarrow \mathcal{T}(\mathbb{R}^d)$ in the translation group of $\mathbb{R}^d$, whose image $\pi(\Gamma)$ has maximal rank $d$. In addition, we require the compatibility condition of $p$ and $\pi$ for each $\gamma \in \Gamma$ and $v\in V$:
\[p(\gamma v)=\pi(\gamma)p(v).\]
\end{definition}

This introduces the static notion of a periodic framework. However, auxeticity requires deformation of the framework. Working in this direction, we want to understand paths in the framework's configuration space.
\begin{definition}
\label{def:oneparameterdeformation}
A one-parameter deformation of a $d$-periodic framework $\mathcal{F}=(G, \Gamma, p, \pi)$ is a smooth family of placements $p_\tau :V\rightarrow \mathbb{R}^d$ parametrized by $\tau \in (-\varepsilon,~\varepsilon)$ for some small $\varepsilon >0$ with $p_0$ defined as the initial placement $p$, satisfying two conditions:
\begin{enumerate}
    \item[a)] The bar-length equations in $b$ remain satisfied and
    \item[b)] Periodicity under $\Gamma$ is maintained via a faithful representation $\pi_\tau:\Gamma\rightarrow \mathcal{T}(\mathbb{R}^d)$. Note that the periodicity lattice $\pi_\tau(\Gamma)$ may change with $\tau$.
\end{enumerate}
After factoring out Euclidean motions, the framework's configuration space denotes the collection of periodic placements in $\mathbb{R}^d$ that satisfy the bar-length equations in $b$. The deformation space of the framework denotes the connected component of the configuration space containing the initial framework $p_0$. 
\end{definition}

Deformation spaces are semi-algebraic sets, implying that notions from algebraic and differential geometry such as singularity, tangent space and dimension apply here as well. The theory of auxetic frameworks -- characterized by perpendicular expansion upon stretching a framework in a chosen direction -- can thus be stated as follows \cite{geometricauxetics}. 
\begin{definition}
\label{def:contraction}
Let $T:\mathbb{R}^d\rightarrow \mathbb{R}^d$ be a linear operator. Define the operator norm of $T$ by
\[||T||=\sup _{|x|\leq 1} |Tx|=\sup _{|x|= 1} |Tx|.\]
$T$ is called a contraction, when $||T||\leq 1$ and a strict contraction, if $||T||<1$.
\end{definition}

Assume now that $\mathcal{F}=(G,\Gamma,p_\tau,\pi_\tau)$ is a one-parameter deformation with $\tau \in (-\varepsilon,\varepsilon)$ of a periodic framework in $\mathbb{R}^d$, as was introduced in \Cref{def:oneparameterdeformation}. Given parameters $\tau_1<\tau_2$, the corresponding one-parameter family of periodicity lattices $\pi_\tau (\Gamma)$ yields a way to compare the framework at $\tau_1$ with the framework at $\tau_2$. Since $\pi_\tau(\Gamma)$ has rank $d$ by assumption, the linear operator $T_{\tau_2\tau_1}$ taking the lattice at $\tau_2$ to $\tau_1$ via
\begin{eqnarray}
\label{eqn:linearoperatorauxetic}
T_{\tau_2 \tau_1} \circ \pi_{\tau_2} = \pi_{\tau_1}
\end{eqnarray}
is unique, leading to the following definition, formalizing auxeticity \cite{geometricauxetics}.

\begin{definition}
\label{def:auxeticpath}
A differentiable one-parameter deformation $(G,\Gamma,p_\tau,\pi_\tau)$ with $\tau \in (-\varepsilon,\varepsilon)$ of a periodic framework in $\mathbb{R}^d$ is an auxetic path, when for any $\tau_1<\tau_2$, the linear operator $T_{\tau_2 \tau_1}$ defined in $(\ref{eqn:linearoperatorauxetic})$ is a contraction.
\end{definition}

\begin{figure}[h!]
    \centering
    \includegraphics[width=1.\linewidth]{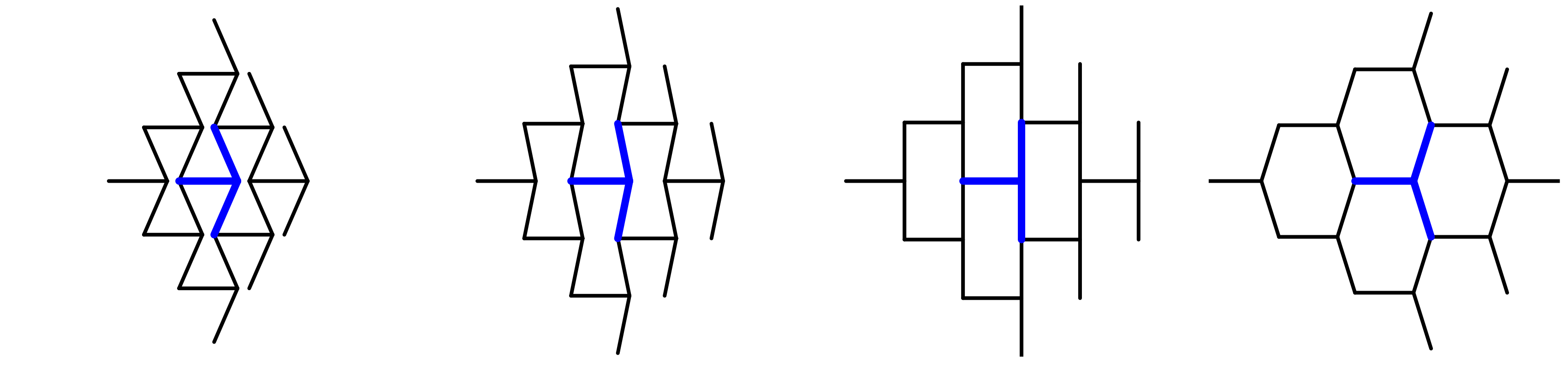}
        \caption{\small All pictures depict 9 unit cells of the honeycomb framework, with one unit cell highlighted in blue. The discussed auxetic behavior becomes clearly visible, when applying a horizontal flex to the framework. At $\tau=0.6$ (l.), the inverted honeycomb is rather compressed and extends over $\tau=0.8$ (c.l.) up to its maximally extended configuration at $\tau=1.0$ (c.r.). After that, the configuration's lateral extension decreases until $\tau=1.3$ (r.) and will continue to do so. In the last picture, the hexagonal honeycomb structure can be observed. }
    \label{fig:auxetichoneycomb}
\end{figure}

To illustrate auxeticity in a periodic structure, we consider a framework based on the honeycomb graph (the vertices and edges of a tiling by regular hexagons) in $\mathbb{R}^2$. This framework exhibits auxetic behavior on a deformation path that is to be defined later. The underlying periodic graph is defined by four vertices $v_1,v_2,v_3,v_4$ and rigid bars $(1,2),~(2,3),~(2,4)$. We define the framework's initial placement as
\[v_1,~v_2,~v_3,~v_4\mapsto (0,0),~(1,0),~(1,1), ~(1,-1).\]
The corresponding periodicity group is generated by $v_3-v_1$ and $v_4-v_1$. Let us first factor out euclidean motions by fixing
\[p(v_1)=(p_{11},~ p_{12})=(0,0) ~\text{ and }~p(v_2)=(p_{21},~ p_{22})=(p_{21},0)\]
for the coordinates $p_{ij}$ of $p(v_i)$. Denoting the edge length of the bar $(i,j)$ by $\ell_{ij}$, the framework's initial placement lets us deduce that $\ell_{ij} = 1$ for all rigid bars, so $p_{21}=1$. This implies that the deformation space is actually a 2-dimensional algebraic set, consisting of all points $(p_{31},~p_{32},~p_{41},~p_{42})$ in $\mathbb{R}^4$ satisfying
\[(p_{31}-p_{21})^2+(p_{32}-p_{22})^2=1 ~\text{ and } ~(p_{41}-p_{21})^2+(p_{42}-p_{22})^2=1.\]
Let us now choose the deformation path parametrized as 
\[\varphi ( \tau ) = \left(\tau, ~\sqrt{2\tau-\tau^2}, ~\tau,~-\sqrt{2\tau-\tau^2}\right)~~~\text{ with }~~~\tau \in (0.5,~1.5).\]
We can quickly check that this path matches \Cref{def:oneparameterdeformation}. The periodicity lattice depends on $\tau$ via 
\[\pi_\tau(\Gamma) = \left(\varphi(\tau)_1,~\varphi(\tau)_2\right) \cdot \mathbb{Z} ~+~ \left(\varphi(\tau)_3,~\varphi(\tau)_4\right)\cdot \mathbb{Z}.\]
According to Equation $(\ref{eqn:linearoperatorauxetic})$, this induces a linear operator
\[T_{\tau_2\tau_1} = \begin{pmatrix}
\frac{\tau_1}{\tau_2} & 0 \\
0 & \sqrt{\frac{2\tau_1-\tau_1^2}{2\tau_2-\tau_2^2}}
\end{pmatrix}.\]
Since $\tau_1,\tau_2\in (0.5,~1.5)$, it holds that $2\tau_i-\tau_i^2>0$ for $i\in \{1,2\}$, implying that $T_{\tau_2\tau_1}$ is well-defined on that interval. We can then check that $||T_{\tau_2\tau_1}||<1$ for $\tau_1,~\tau_2 \in (0.5,1)$ and $||T_{\tau_2\tau_1}||>1$ for $\tau_1,\tau_2 \in (1,1.5)$. Consequently, $\varphi(\tau)$ is an auxetic deformation path for $\tau \in (0.5,1]$ according to \Cref{def:auxeticpath}, with selected deformations illustrated in \Cref{fig:auxetichoneycomb}.

Now, we are ready to introduce the concept of a \emph{tensegrity framework}. Given an undirected graph $G=(V,E)$ with labeled vertices $V$ and edges $E=\mathcal{B}\sqcup \mathcal{C} \subseteq \binom{V} {2}$, 
we call the edges $ij \in \mathcal{B}$ \emph{rigid bars} of length $\ell_{ij}$ and the edges $ij \in \mathcal{C}$ \emph{elastic cables} with corresponding natural resting length $r_{ij}$ and constant of elasticity $c_{ij}$. 
A \emph{tensegrity framework} is a graph $G=(V,E; \mathcal{B},\mathcal{C})$ together with an embedding $p:V\rightarrow \mathbb{R}^d$, assigning a point in $\mathbb{R}^d$ to each vertex of $G$ \cite{catastrophetensegrity, TensegrityFrameworksOriginal}. We call $p$ a \emph{realization} of $G$ and denote the coordinates of the $n$ nodes by $p_1=(p_{11},\dots ,~p_{1d}),~\dots,~p_n$. For every $ij \in \mathcal{B}$, we assign the bar constraint polynomial
\[b_{ij}:= \sum_{k=1}^d (p_{ik}-p_{jk})^2-\ell_{ij}^2\]
and denote by $b$ the polynomial system consisting of the $b_{ij}$ for $ij\in \mathcal{B}$. The cables can be modeled as one-sided Hookean springs, giving rise to the cables' potential energy

\[q_{ij}:=\frac{c_{ij}}{2}\left(\max \{0,~ ||p(i)-p(j)||-r_{ij}\}\right)^2, ~~~Q = \sum_{ij\in\mathcal{C}}q_{ij}.\]
The variables $p_{ik},~\ell_{ij},~r_{ij},~c_{ij}$ either lie in the set of externally determined \emph{control parameters} $Y=\mathbb{R}^{m_1}$, or in the \emph{internal variables} $X=\mathbb{R}^{m_2}$. 
Some of $Y$'s elements are fixed, while we let others vary in some subset $\Omega\subset Y$. We call a configuration $(x,y)\in X\times Y$ \emph{stable} or an \emph{equilibrium}, if $x$ is a strict local minimum of the energy function $Q$ with respect to the algebraic set $\mathcal{V}(b)$. 
Consequently, for fixed $y\in Y$ the search for equilibrium configurations in this tensegrity framework can be modeled as a constraint optimization problem

\begin{eqnarray*}
&~&\min_{x \in X} Q(x,y)\\
&\text{s.t.}& ~b_{ij}(x,y)=0 ~~~~\forall ~ ij\in \mathcal{B}.
\end{eqnarray*}

Naturally, concepts like periodicity and auxeticity that were previously discussed can be transferred to this setting: While periodicity can be extended by considering the two subgraphs $(V, \mathcal{B})$ and $(V,\mathcal{C})$ with identical periodicity group $\Gamma$ (Def. \ref{def:periodicgraph}), we expect a deformation path, parametrized by the varying control parameters $\Omega \subset Y$, to be in equilibrium at each step (Def. \ref{def:oneparameterdeformation}). Similarly, an auxetic deformation path is described by the contraction $T_{\tau_2\tau_1}$ (Def. \ref{def:auxeticpath}).

\section{Negative Poisson's ratio}
\label{section:previousmodel}

From a materials science perspective, the Poisson's ratio is a mechanical property that describes a material's deformation behavior under loading. It is defined as the ratio of the lateral strain (change in width or thickness) to the axial strain (change in length) of a material when subjected to an applied load. It is an infinitesimal quantity of a linear elastic material. 

The tensegrity structures that we are analyzing here have non-linear behavior exhibited over large deformations. We refer then to the \emph{instantaneous Poisson's ratio}, which allows us to consider behavior over a larger range of deformations. For a given set of discrete timesteps, $t_0 < t_1 < \dots < t_N$, it is defined in terms of the engineer's strain \cite{poissonratio}:
\begin{eqnarray}
\label{def:poissonratio}
\nu_{xy} &=& -\frac{e_y}{e_x}\\[1mm]
\text{with}~~~~~e_x &=& \frac{(L_x)_{t_j} - (L_x)_{t_{j-1}}}{(L_x)_{t_{j-1}}}~~~~~\text{and}~~~~~e_y = \frac{(L_y)_{t_j} - (L_y)_{t_{j-1}}}{(L_y)_{t_{j-1}}}.\nonumber
\end{eqnarray}
Here, $x$ is the direction of applied strain and $y$ is an orthogonal direction. $(L_x)_t$ and $(L_y)_t$ are the lattice extensions in the $x$- and $y$-directions taken at timestep $t$. Materials with negative Poisson's ratio are called \emph{auxetic}.

Notice that we have now introduced two different concepts of auxeticity -- the geometric \Cref{def:auxeticpath} and the instantaneous Poisson's ratio over large deformations in Equation (\ref{def:poissonratio}) --, which we wish to compare. The geometric definition of auxeticity considers shearing as well as orthogonal extensions of the unit cell, in contrast to the instantaneous Poisson's ratio shown above, whose increase in all orthogonal directions is not sufficient for a deformation path to be auxetic geometrically \cite{geometricauxetics}. Nevertheless, the instantaneous Poisson's ratio is easier to compute in certain experimental and numerical contexts, so often times it is preferred in practice. The following proposition relates these two concepts. 
\begin{proposition}
\label{prop:auxeticityandpoissonratio}
Assume that $G = (V,E,p,\pi)$ is a $d$-periodic tensegrity framework. If the deformation path $(p_\tau,\pi_\tau)$ induced by stretching the framework in a fixed direction is auxetic in the sense of \Cref{def:auxeticpath}, then the Poisson's ratio is non-positive for any discretizations $-\varepsilon < t_0 < \dots < t_N < \varepsilon$. Conversely, if the lattice $\pi_\tau$ has an orthogonal unit cell for all $\tau$ and the Poisson's ratio is non-positive for all discretizations of $(-\varepsilon,\varepsilon)$, then the deformation path is auxetic.
\begin{proof}
Let $\Lambda_\tau$ denote the matrix of generators for the periodicity lattice $\pi_\tau(\Gamma)$ corresponding to the framework's deformation path. After factoring out Euclidean motions, applying the assumption that the framework is stretched in a fixed direction and potentially reordering, the matrix of generators becomes
\begin{eqnarray}
\label{eqn:latticeexpression}
\Lambda_\tau = \begin{pmatrix} \tau + \ell_0 & 0 & \dots & 0\\
c_{12}(\tau) & c_{22}(\tau) & ~& \vdots\\
\vdots & ~ & \ddots & 0\\
c_{1d}(\tau) & \dots & c_{(d-1)d}(\tau) & c_{dd}(\tau)\end{pmatrix}
\end{eqnarray}
for continuous functions $c_{ij}$ -- a topological path is continuous -- and the initial extension $\ell_0 \neq 0$. Since $\Lambda_\tau$ is lower-triangular, for any $\tau_1<\tau_2$, the linear operator $T_{\tau_2\tau_1}$ defined in (\ref{eqn:linearoperatorauxetic}) taking the lattice $\pi_\tau$ from time $\tau_2$ to $\tau_1$ has lower-triangular matrix representation, too. Therefore, its eigenvalues are on the diagonal \cite[p.152]{axlerlinearalgebra} and are precisely given by 
\[\lambda(T_{\tau_2\tau_1}) = \left(\frac{\tau_1+\ell_0}{\tau_2+\ell_0},~\frac{c_{22}(\tau_1)}{c_{22}(\tau_2)},~\dots,~\frac{c_{dd}(\tau_1)}{c_{dd}(\tau_2)}\right).\] 
By the Spectral Radius Theorem \cite[p.347]{matrixanalysis}, the eigenvalues' absolute values are bounded above by any matrix norm. Therefore, assuming that the deformation path is auxetic immediately proves that $|c_{ii}(\tau_2)|\geq|c_{ii}(\tau_1)|$ for $\tau_2>\tau_1$. By the framework's $d$-periodicity, $\Lambda_\tau$ has full rank for any $\tau$. Therefore, $c_{ii}(\tau)\neq 0$, so $e_i \geq 0$ for each $i$, implying that the instantaneous Poisson's ratio is non-positive for any discretization of the interval $(-\varepsilon, \varepsilon)$.

Conversely, assuming that the Poisson's ratio is non-positive, we can deduce that $\ell_0>\varepsilon$, by choosing a positive branch and by the framework's $d$-periodicity. Therefore, for $\tau_2> \tau_1$

\begin{eqnarray*}
0 &\geq& \nu_{1i} = -\frac{e_i}{e_1} = -\frac{\tau_1+\ell_0}{\lvert c_{ii}(\tau_1)\lvert} \cdot \frac{\lvert c_{ii}(\tau_2)\lvert-\lvert c_{ii}(\tau_1)\lvert}{\tau_2-\tau_1}~~~~\forall i>1,~~~~~~~~~\text{so}\\
0 &\leq& \frac{\lvert c_{ii}(\tau_2)\lvert-\lvert c_{ii}(\tau_1)\lvert}{\lvert c_{ii}(\tau_1)\lvert} ~~~~\forall i>1.
\end{eqnarray*}

In particular, $|c_{ii}(\tau_2)|\geq |c_{ii}(\tau_1)|$. By assumption, the unit cell is orthogonal, so after potentially reordering, the lattice $\Lambda_\tau$ can be chosen as a diagonal matrix at each step. Together with the framework's $d$-periodicity, this implies that $c_{ii}(\tau)\neq 0$ for each $\tau$. By the deformation path's continuity, we find that $c_{ii}(\tau_1)$ and $c_{ii}(\tau_2)$ have the same sign. Consequently, $c_{ii}(\tau_1)/c_{ii}(\tau_2)\leq 1$.
Therefore, all eigenvalues of $T_{\tau_2\tau_1}$ are smaller than $1$, yet still positive. Previous observations then imply that the linear operator $T_{\tau_2\tau_1}$ is also diagonal. As a diagonal matrix is normal, the absolute values of its eigenvalues and singular values agree. Thus, the operator norm induced by the Euclidean norm is equal to the largest eigenvalue of $T_{\tau_2\tau_1}$, which is at most 1, making the linear operator a contraction. As a result, the deformation path is auxetic by \Cref{def:auxeticpath}.
\end{proof}
\end{proposition} 

However, \Cref{prop:auxeticityandpoissonratio} is not exactly useful in determining, whether an arbitrary framework is auxetic. The assumption that the lattice stays orthogonal at every time step often is too much to ask for. For this reason, it seems beneficial to give another sufficient criterion for auxeticity in terms of the lattice generators' coordinates -- at least on a given discretization. Basically, we formulate a bound on the size of the lattice's off-diagonal entries relative to the diagonal entries, dependent on the step size with which the deformation path is discretized. As our setting is motivated by real-world examples, we only prove the result for $3$-periodic frameworks, though the proposition can be analogously proven for $2$-periodic structures.

\begin{proposition}
\label{prop:inequalitiesforauxeticity}
Assume that $G = (V,E,p,\pi)$ is a $3$-periodic tensegrity framework. Let the deformation path $(p_\tau,\pi_\tau)$ with $\tau \in (-\varepsilon, \varepsilon)$ be induced by stretching the framework in a fixed direction along a given discretization $-\varepsilon < t_0 < \dots < t_N < \varepsilon$. Further assume that the absolute values of  $\Lambda_\tau$'s (cf. (\ref{eqn:latticeexpression})) off-diagonal entries $|c_{jk}(\tau)|$ are bounded above by $\alpha|c_{ii}(\tau)|$ for some $\alpha>0$, all grid points $\tau$ and $i\in \{1,2,3\}$. Finally, assume that 
\begin{eqnarray*}
\left(\frac{c_{11}(\tau_1)}{c_{11}(\tau_2)}\right)^2 & \leq & \frac{1-3\alpha-11\alpha^2-12\alpha^3-4\alpha^4}{1+2\alpha^2},\\[1mm]
\left(\frac{c_{22}(\tau_1)}{c_{22}(\tau_2)}\right)^2 & \leq & \frac{1-3\alpha-7\alpha^2-3\alpha^3}{1+\alpha+\alpha^2+\alpha^3}~~~~~\text{ and }\\[1mm]
\left(\frac{c_{33}(\tau_1)}{c_{33}(\tau_2)}\right)^2 & \leq & \frac{1-2\alpha-\alpha^2}{1+2\alpha+\alpha^2}
\end{eqnarray*} 
for grid points $\tau_2>\tau_1$. Then, the linear operator $T_{\tau_2\tau_1}$ defined in Equation (\ref{eqn:linearoperatorauxetic}) is a contraction for all grid points $\tau_2>\tau_1$.
\begin{proof}
    Analogous to \Cref{prop:auxeticityandpoissonratio}'s proof (cf. Equation (\ref{eqn:latticeexpression})), after removing Euclidean motions the linear operator taking the lattice $\pi_\tau$ from time $\tau_2$ to $\tau_1$ can be expressed as the lower-triangular matrix
    \[T_{\tau_2\tau_1} = \begin{pmatrix} 
    a_1 & 0 & 0 \\
    a_2 & a_4 & 0\\
    a_3 & a_5 & a_6\end{pmatrix},\]
with 
\begin{eqnarray*}
a_i &=& \frac{c_{ii}(\tau_1)}{c_{ii}(\tau_2)} ~~~\text{ for }~~~i\in \{1,4,6\},\\
a_{3i-1} &=& \frac{c_{i,i+1}(\tau_1)c_{i+1,i+1}(\tau_2)-c_{i,i+1}(\tau_2)c_{i+1,i+1}(\tau_1)}{c_{i,i}(\tau_2)c_{i+1,i+1}(\tau_2)} ~~~\text{ for }~~~i\in \{1,2\} ~~ \text{ and}\\
a_3 &=& \frac{c_{13}(\tau_1)c_{33}(\tau_2)-a_5c_{12}(\tau_2)c_{33}(\tau_2)-c_{13}(\tau_2)c_{33}(\tau_1)}{c_{11}(\tau_2)c_{33}(\tau_2)}.
\end{eqnarray*}
To show that $T_{\tau_2\tau_1}$ is a contraction, we need to show that the largest eigenvalue of $T_{\tau_2\tau_1}^TT_{\tau_2\tau_1}$ is at most $1$. To calculate the matrix' eigenvalues, we want to utilize the Greshgorin Circle Theorem \cite[Satz II]{gerschgorin}, which states that the eigenvalues of $T_{\tau_2\tau_1}^TT_{\tau_2\tau_1}$ lie in (complex) circles whose centers are the matrix' diagonal entries and whose radii are given by the sum of the off-diagonal row entries' absolute values. With the claim's assumptions and the triangle inequality, we can now compute the following:
\begin{eqnarray*}
    |a_{3i-1}| &\overset{\Delta}{\leq}& \left\lvert \frac{c_{i,i+1}(\tau_1)}{c_{i,i}(\tau_2)}\right\lvert + \left\lvert \frac{c_{i,i+1}(\tau_2)\cdot c_{i+1,i+1}(\tau_1)}{c_{i,i}(\tau_2)\cdot c_{i+1,i+1}(\tau_2)}\right\lvert  \\
    &\leq& \alpha \cdot \left( \left\lvert a_{3i-2} \right\lvert + \left\lvert a_{2i+2} \right\lvert\right) ~~~~~\text{ for } i \in \{1,2\} ~~~~~\text{ and }\\[3mm]
    |a_3| &\overset{\Delta}{\leq} & \left\lvert \frac{c_{13}(\tau_1)}{c_{11}(\tau_2)}\right\lvert + |a_5|\cdot \left\lvert \frac{c_{12}(\tau_2)}{c_{11}(\tau_2)}\right\lvert + |a_6|\cdot \left\lvert \frac{c_{12}(\tau_2)}{c_{11}(\tau_2)}\right\lvert \\
    &\leq& \alpha \cdot\left( |a_1| + |a_6| + \alpha \cdot \left(|a_4|+|a_6|\right)\right).
\end{eqnarray*}
These results finally enable us to calculate the row sums in absolute values:
\begin{eqnarray*}
\text{I}:&~& \left({a_1}^2+{a_2}^2+{a_3}^2\right) + |a_2a_4+a_3a_5|+|a_3a_6|~ \\
&~&~~~\leq ~|a_4|^2\cdot\left(\alpha+\alpha^2+\alpha^3+\alpha^4\right)+|a_1|^2\cdot\left(1+2\alpha^2\right)+|a_4a_6|\cdot\left(2\alpha^2+4\alpha^3+2\alpha^4\right)\\
&~&~~~+|a_1a_4|\cdot\left(\alpha+3\alpha^2+2\alpha^3\right) +|a_1a_6|\cdot\left( \alpha+3\alpha^2+2\alpha^3 \right)+ |a_6|^2\cdot\left(2\alpha^2+3\alpha^3+\alpha^4\right)  \\[1mm]
\text{II}:&~& \left({a_4}^2+{a_5}^2\right) + |a_2a_4+a_3a_5|+|a_5a_6|~\leq ~|a_4|^2\cdot \left( 1+\alpha+\alpha^2+\alpha^3\right) \\
&~&~~~+ |a_1a_4|\cdot\left(\alpha+\alpha^2\right) + |a_1a_6|\cdot \alpha^2 + |a_4a_6|\cdot \left(\alpha+3\alpha^2+2\alpha^3\right) + |a_6|^2\cdot \alpha \cdot (1+\alpha)^2\\[1mm]
\text{III}:&~&\left({a_6}^2\right)+|a_3a_6|+|a_5a_6|~\leq ~ |a_6|\cdot\left(\alpha ~ |a_1| + (\alpha+\alpha^2)~|a_4| + (1+2\alpha+\alpha^2)~|a_6|\right)
\end{eqnarray*}
By the assumptions on $(a_1)^2,~(a_4)^2$ and $(a_6)^2$, the values of I$-$III are bounded above by $1$, so Greshgorin's Circle Theorem implies that the eigenvalues of $T_{\tau_2\tau_1}^TT_{\tau_2\tau_1}$ are bounded above by $1$. \Cref{def:contraction} then implies that the linear operator $T_{\tau_2\tau_1}$ is a contraction.
\end{proof}
\end{proposition}

\begin{remark}
\Cref{prop:inequalitiesforauxeticity} immediately begs the question, whether the assumptions on $c_{ii}(\tau_1)/c_{ii}(\tau_2)$ are reasonable. Indeed, for $0\leq \alpha<0.183$, the expressions on the right of the inequalities are positive. Nevertheless, having $\alpha$ close to $0.183$ would mean that $c_{ii}(\tau_2)$ is insanely large relative to $c_{ii}(\tau_1)$. More sensibly, by choosing a constant step length of $10^{-2}$, maximum extension $c_{ii}(0)=1$ for $i\in \{1,2,3\}$ and constant Poisson's ratio $-1$ we can deduce that $\alpha \leq 4\cdot 10^{-3}$. In particular, the smaller the step size, the smaller $\alpha$ will be.
\end{remark}

\section{Constructing a stable periodic tensegrity}
\label{section:robustmodel}
\begin{figure}[h!]
\centering
\begin{minipage}{0.43\textwidth}
    \includegraphics[width=\textwidth]{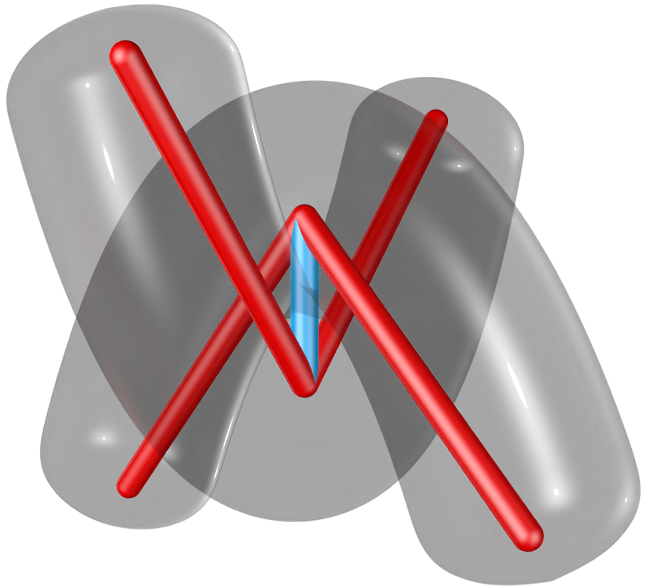}
    \label{fig:tensegritymodel}
\end{minipage}
\hfill
\begin{minipage}{0.43\textwidth}
    \includegraphics[width=\textwidth]{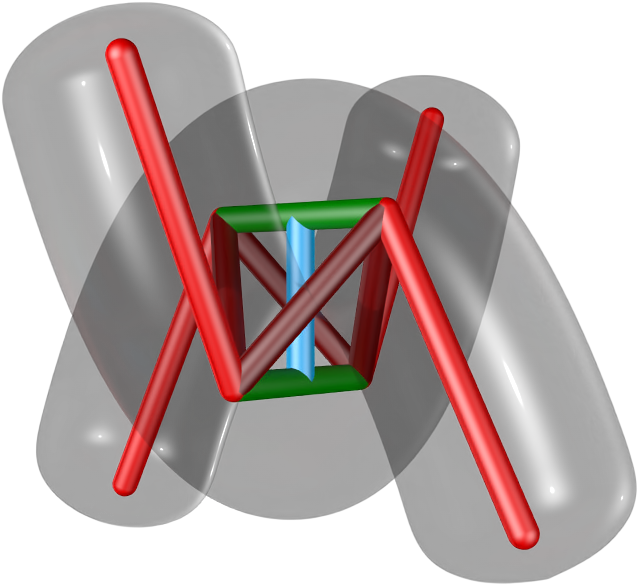}
    \label{fig:tetrahedralmodel}
\end{minipage}
\caption{\small The original tensegrity model (l.) and the tetrahedral tensegrity model (r.) for two filaments in tight, orthogonal contact.}
\label{fig:claspfigures}
\end{figure}

In a previous study, the $\Pi^+$ rod packing \cite{Evans2011} (\Cref{fig:bmnsgn}) was transformed into a tensegrity by adding balancing equations \cite{reentranttensegrity}. In that model, each contact between filaments is described by an incompressible bar of radius $r>0$ normal to the contact, connecting the central axes of the cylinders. As depicted in \Cref{fig:claspfigures}l., flexible elastic cables are placed along the cylinders' central axes to model the filaments' elasticity. The tensegrity structure was shown to exhibit auxetic behavior, albeit with very unstable numerics for the initial phase of the deformation parametrized by stretching the corresponding framework along the $x$-axis. We suspect that the unstable behavior is induced by the instability of the contact bar between two filaments. 

Theoretical studies and recent experiments of the contact set of two orthogonal tubes in a clasp configuration show that the system balances through a more complicated intersection than just a single point \cite{criticalitygehring,clasptightorthogonalcontact, starostin}: Two filaments in tight contact meet in a contact surface, on which the pressure varies. This idea suggest that the contact between two filaments could be better modeled as a tetrahedron of constraints, encapsulating the four points of highest pressure. We describe this model here.

To build the new model for two filaments in contact, denote by $p_{11}, p_{12}, p_{21}, p_{22} \in \mathbb{R}^3$ four points in space. If these four points $p_{ij}$ are in general position, they form a tetrahedron with $12$-dimensional configuration space. Assume now that $p_{i1}$ and $p_{i2}$ lie in the center of the same filament. 
We place cables $c_{ij}$ along the cylinders' central axes,  following the same notation as the point they are attached to, so $c_{12}$ would be attached to $p_{12}$. Assume that the cylinders' radius is equal to $r > 0$. Analogous to the previous model, we want the cylinders' centers to be $2r$ apart, meaning that we place a bar of length $2r$ in the tetrahedron's center, depicted in blue in \Cref{fig:claspfigures}r. This gives rise to the first constraint,
\begin{eqnarray} 
\label{eqn:centerbarlength}
\left\lvert \frac{p_{11} + p_{12}}{2} - \frac{p_{21} + p_{22}}{2} \right\lvert = 2 r.
\end{eqnarray}
Furthermore, this constraint should not only keep the cylinders apart at one position, but rather throughout the entire contact. Because of this and for symmetry reasons, we require that the center bar (blue) is orthogonal to both of the tetrahedron's sides (green) by setting
\begin{eqnarray}
\label{eqn:orthogonalcenterbar}
 \left\langle\frac{p_{11} + p_{12}}{2} - \frac{p_{21} + p_{22}}{2},~~p_{i1} - p_{i2} \right \rangle &=& 0 ~~~~~~~~\text{for } i \in \{1,2\}.
\end{eqnarray}
To maintain its symmetry, we assume that the angles between the incoming cables $c_{i1}, c_{i2}$ and the tetrahedron's sides $p_{i2}-p_{i1}$ (green) are consistent. This translates to the equation
\begin{eqnarray}
\label{eqn:sameorientation}
\frac{\left\langle\frac{p_{11} + p_{12}}{2} - \frac{p_{21} + p_{22}}{2},~c_{i1}\right \rangle}{\left\lvert c_{i1} \right\lvert } = \frac{\left\langle\frac{p_{11} + p_{12}}{2} - \frac{p_{21} + p_{22}}{2},~c_{i2}\right \rangle}{\left\lvert c_{i2} \right\lvert }~~~~~~~~\text{for } i \in \{1,2\}. 
\end{eqnarray}
Additionally, we require that the tetrahedron does not twist in unforeseen ways by assuming that the incoming cables $c_{i1}$ and $c_{i2}$ lie in the same plane as the tetrahedron's sides (green). This can be expressed by the equation
\begin{eqnarray}
\label{eqn:orthogonalcables}
\left\langle c_{i1}\times c_{i2}, ~ p_{i1}-p_{i2} \right \rangle = 0 ~~~~~~~~\text{for } i \in \{1,2\}. 
\end{eqnarray}
Lastly, the tetrahedron can open. For that reason, the tetrahedron's sides (green) are not fixed, but should rather depend on the attached cables. Denote the angle between the cables $c_{i1}$ and $c_{i2}$ by $\theta_i \in (0,\pi)$. We can express the varbiable bar's length (green) as $d_i = 2 \sin \left(\theta_i/2\right)$ (cf. \cite[Thm. 9.5]{criticalitygehring}). By applying trigonometric identities, we then calculate
\[d_i = 2 \cdot \sin \left(\frac{\theta_i}{2}\right) = 2 \cdot \text{sgn}\left(\sin \left(\frac{\theta_i}{2}\right)\right)\cdot \sqrt{\frac{1-cos(\theta_i)}{2}} = 2 \cdot \text{sgn}\left(\sin \left(\frac{\theta_i}{2}\right)\right)\cdot \sqrt{\frac{1-\frac{\langle c_{i1},~c_{i2}\rangle}{|c_{i1}|\cdot |c_{i2}|}}{2}}\]
To be meaningful, this expression needs to be scaled with respect to the tetrahedron's size, so we multiply the right hand side of the above equation by $r$. As $\theta_i \in (0,\pi)$, it holds that $\sin(\theta_i/2)\geq 0$, enabling us to rewrite the expression in point coordinates, i.e. $d_i=|p_{i1}-p_{i2}|$:
\begin{eqnarray}
\label{eqn:variablebarlength}
\langle c_{i1},~c_{i2}\rangle  ~ = |c_{i1}|\cdot |c_{i2}|\cdot \left(1-\frac{|p_{i1}-p_{i2}|^2}{2\cdot r^2} \right).
\end{eqnarray}
Besides the four incoming cables $c_{ij}$ (red), we also place four internal cables (dark red) on the remaining edges of the tetrahedron to penalize twisting. In total, for generic placements $p$, there are $4\cdot 3 - 9 = 3$ degrees of freedom per tetrahedron. They can be parametrized by the lengths of the two variable bars (green) and the relative orientation of these two bars. The constraints are smooth almost everywhere, except when the external cables or variable bars have zero length, making Riemannian optimization algorithms \cite[p.62]{boumal2020intromanifolds} feasible.

To formalize the tetrahedra and relate them to the geometric theory introduced in \Cref{section:geometricauxetics}, we now consider an embedded hypergraph $G_\mathcal{T}=(V\sqcup V,\mathcal{T},p,\pi)$ with the set of ordered hyperedges defined by the tetrahedron's edges. We write $i_1$ for an element of the first $V$ in the disjoint union, $i_2$ for an element of the second $V$ in the disjoint union and $i$ for an element of the original set of vertices $V$. Recall that in the original model, each vertex has two cables attached and one bar.
\[\mathcal{T} = \{\{(h_2,i_1,i_2,k_1),(m_2,j_1,j_2,n_1)\} : ~\text{for }ij\in \mathcal{B} \text{ with } hi,ik, mj,jn\in \mathcal{C}\}\]
For an element $T\in \mathcal{T}$ we write $T_{i,j}$ for the $j$-th element in the $i$-th entry of $T$ ($i\in \{1,2\}$, $j\in\{1,2,3,4\}$) and $p_{ij}$ for $p(T_{i,j})$. With the notation in place, we can now summarize the objective function derived from the cables envisioned as one-sided Hookean springs and the constraints (\ref{eqn:centerbarlength}) - (\ref{eqn:variablebarlength}) in the nonlinear optimization problem
\begin{eqnarray}
\label{eqn:optimizationproblemnew}
&\underset{{p:V\rightarrow\mathbb{R}}}{\min}& \sum _{T\in \mathcal{T}} \left( \sum_{\substack{i\in\{1,2\},\\j\in\{1,3\}}}\frac{c_{ij}}{2}\left(\max \{0,~ ||p_{ij}-p_{ij+1}||-r_{ij}\}\right)^2 + \right.\\
&~&\quad\quad~\left.\sum_{i,j\in\{2,3\}}\frac{c_{ij}}{2}\left(\max \{0,~ ||p_{1i}-p_{2j}||-r_{ij}\}\right)^2 \right) \nonumber\\[2mm]
&\text{s.t.}& \left\lvert \frac{p_{12} + p_{13}}{2} - \frac{p_{22} + p_{23}}{2} \right\lvert = 2 r,\nonumber\\
&~& \left\langle \frac{p_{12} + p_{13}}{2} - \frac{p_{22} + p_{23}}{2},~p_{i2} - p_{i3} \right \rangle = 0,\nonumber\\
&~& \frac{\left\langle\frac{p_{12} + p_{13}}{2} - \frac{p_{22} + p_{23}}{2},~p_{i2}-p_{i1}\right \rangle}{\left\lvert p_{i2}-p_{i1} \right\lvert } = \frac{\left\langle\frac{p_{12} + p_{13}}{2} - \frac{p_{22} + p_{23}}{2},~p_{i4}-p_{i3}\right \rangle}{\left\lvert p_{i4}-p_{i3} \right\lvert }, \nonumber\\
&~& \left\langle (p_{i2}-p_{i1})\times (p_{i4}-p_{i3}), ~ p_{i3}-p_{i2} \right \rangle = 0~~~~~~\text{ and}\nonumber\\
&~& \langle p_{i2}-p_{i1},~p_{i4}-p_{i3}\rangle  ~ = |p_{i2}-p_{i1}|\cdot |p_{i4}-p_{i3}|\cdot \left(1-\frac{|p_{i3}-p_{i2}|^2}{2\cdot r^2} \right)
\nonumber
\end{eqnarray}
with constraints taken for each $T\in \mathcal{T}$ and $i \in \{1,2\}$. For the incoming cables, we choose the resting length $r_{ij}=0.1$ and for the cables inside the tetrahedron $r_{ij}=2r$. Analogously, the cables' constant of elasticity is $c_{ij}=1$ for incoming cables and $c_{ij}=30$ for the interior cables to strongly penalize deviations from the orthogonal case. 
While the left sum in the objective function corresponds to the incoming cables $c_{ij}$, the right part corresponds to the cables on the tetrahedron's sides. The constraints are obtained by inserting the representation from $\mathcal{T}$ into Equations $(\ref{eqn:centerbarlength})-(\ref{eqn:variablebarlength})$.

This outlines the set up of a tetrahedron of constraints, providing a stable tensegrity structure. The challenge with this model is to generate a set of constraints that make the tetrahedron robust with respect to perturbations, while it remains flexible enough to adjust depending on the incoming cables to accommodate for the different ways two filaments can be in contact. Our use of the orthogonal clasp's contact set in the tensegrity structure's design means that our model best approximates those structures with orthogonal cylinders. Some twisting of the tetrahedron is allowed (but penalized), which extends beyond these constraints, though a more comprehensive study of non-orthogonal cylinders is left for future research.

\section{Riemannian Optimization}
\label{section:riemannianoptimization}
With the well-defined nonlinear optimization problem (\ref{eqn:optimizationproblemnew}) for our tensegrity in place, we utilize a robust optimization algorithm to solve it. A natural approach to finding such an algorithm is using Lagrange multipliers \cite[p.320ff.]{nocedal2006}, as we are exclusively dealing with equality constraints. This leads to a square polynomial system that can be solved using classical methods from numerical algebraic geometry, such as Homotopy Continuation \cite{Wampler2013bertini, HC.jl, Hom4PSArticle}. This method can be used to find \emph{all} isolated solutions of a polynomial system with theoretical guarantees. However, Bernstein's Theorem predicts that the number of solutions can grow exponentially in the amount of variables \cite{coxlittleoshea}. Since the system involved in these tensegrity structures is exceeding 100 equations, the global approach is infeasible.

For that reason, instead of trying to find all solutions, we only intend to find one solution. Indeed, finding one configuration in equilibrium suffices to model the tensegrity's behavior. There is a myriad of algorithms to find an optimum, e.g. conjugate gradient or quasi-Newton \cite[p.101ff.]{nocedal2006}, interior point \cite[p.392f.]{nocedal2006} and augmented Lagrangian \cite[p.497f.]{nocedal2006}. However, such methods usually require a good initialization, which we were not able to provide in our experiments. Consequently, we could not get them to converge to a critical point. 

The methods' main drawback for our setting is that they try to balance the constraints with the objective function, while we really want the constraints to be satisfied at every iteration to prevent catastrophic behavior in the framework (e.g. \cite{catastrophetensegrity}). As mentioned in \Cref{section:robustmodel}, the constraints are smooth almost everywhere, so Riemannian optimization techniques become feasible. The general idea is to consider the constraint set as an implicitly defined manifold and iteratively apply the exponential map to descent directions to stay on the manifold. However, as we do not have access to the exponential map of general manifolds, it is approximated with a retraction \cite[p.46]{boumal2020intromanifolds}. Still, retraction maps are not readily available for general manifolds. Nevertheless, as the manifold is embedded in $\mathbb{R}^d$, the closest point on the manifold with respect to the Euclidean distance is actually a retraction \cite{AbsMal2012}. As the problem of finding the closest point on an algebraic variety can be expressed as a square polynomial system with finitely many solutions \cite{euclideandistancedegree}, we can even use the previously mentioned homotopy continuation methods to solve it.

At first, this may seem absurd: We begin with a constraint optimization problem on the manifold, to minimize an objective function. Then we suggest iteratively solving another constraint optimization problem on the manifold with the goal of finding the closest point in each step. In particular, in general there is no explicit formula. Fortunately, we know a solution to the closest point problem at any point on the manifold, namely the point itself. Then, we can deform this solution via homotopy to obtain the next point, the retraction \cite{euclideandistanceretraction}. The homotopy, in turn, is computed by a predictor-corrector scheme. While state-of-the-art algorithms use higher-order Runge-Kutta methods as predictor \cite{Wampler2013bertini, HC.jl, Hom4PSArticle}, for the sake of brevity we describe a simpler scheme consisting of an Euler step to guess a solution at the next parameter and Newton's method to correct the step back to the manifold. This path-tracking scheme is summarized by \Cref{alg:path-tracking}.

What remains is choosing suitable step sizes and descent directions. For the sake of simplicity, we choose steepest descent projected to the tangent space \cite[p.20f.]{nocedal2006} as step direction and backtracking line search \cite[p.37]{nocedal2006} to determine the next step size. Finally, this enables the use of the Riemannian gradient descent algorithm \cite[p.62f.]{boumal2020intromanifolds}. Conveniently, Zoutendijk's Theorem \cite{Ring2012} guarantess that this algorithm converges. The package \texttt{HomotopyOpt.jl}\footnote{\url{https://github.com/matthiashimmelmann/HomotopyOpt.jl}} is described in greater detail in \cite{euclideandistanceretraction} and implements the Riemanian gradient descent algorithm with Euclidean distance retraction. We use it to find an energy minimum of the nonlinear optimization problem (\ref{eqn:optimizationproblemnew}). In particular, this procedure yields a critical point of the Lagrange multiplier function corresponding to (\ref{eqn:optimizationproblemnew}) and equivalently, a zero of the Lagrange multiplier system \cite[p.321]{nocedal2006}. 

Subsequently, we want to investigate the framework's displacement that occurs when expanding it in a fixed direction. As this extension can be realized as a one-parameter deformation path (cf. Def. \ref{def:oneparameterdeformation}), we can view it as a homotopy of polynomial systems. To solve it, we can again employ the path-tracking \Cref{alg:path-tracking}, finally revealing the optimization pipeline described in pseudocode by \Cref{alg:rodpacking}, that is used to find deformation paths of tensegrity frameworks.

 \begin{algorithm}[H]
 \caption{Path-tracking Algorithm (cf. \cite{sommesewampler})}
 \label{alg:path-tracking}
 \begin{algorithmic}[1]
\STATE{\textbf{Input: }A homotopy $H(x,t)$ with $x(1)$ known, a grid $1=t_0>t_1>\dots>t_N=0$.}
\STATE{\textbf{Output: }A solution of $H(x,0)=0$.}
\STATE{Set $w_0=x(1)$.}
\FOR{$j \in \{0,\dots,N-1\}$}
\STATE{Predict $w=w_j+\Delta x$ by solving $\frac{\partial H}{\partial x}(w_j,t_j)\Delta x = -\frac{\partial H}{\partial t}(w_j,t_j) h$ (Euler's method).}
\STATE{Correct $w$ by solving $H(x,t_{j+1})=0$ (Newton's method), yielding $w_{j+1}$.}
\ENDFOR
\RETURN{$w_N=x(0)$}
 \end{algorithmic}
\end{algorithm}
 \begin{algorithm}[H]
 \caption{Tensegrity deformation }
 \label{alg:rodpacking}
 \begin{algorithmic}[1]
 \STATE{\textbf{Input: }A rod packing's contact graph (cf. \Cref{section:previousmodel}), amount of steps $N$ and a maximum extension $T>0$. }
\STATE{\textbf{Output: }An animation of the one-parameter deformation of the tensegrity framework induced by stretching the rod packing in $x$-direction.}
\item[]
\STATE{Replace each rod with a tetrahedron in accordance with the nonlinear optimization problem (\ref{eqn:optimizationproblemnew}), yielding an objective function $Q$ and a set of constraints $g$.}
\STATE{Apply the Riemannian gradient descent algorithm with Euclidean distance retraction to find an energy minimum at $x_0$.}
\STATE{Formulate the Lagrange multiplier function $\mathcal{L}(x,\lambda;\tau) = Q(x;\tau)+\lambda^Tg(x;\tau)$, yielding a deformation path parametrized by $\tau\in [0,T]$ with $\nabla_{x,\lambda} \mathcal{L}(x_0, \lambda_0;0)=0$.}
\STATE{Apply \Cref{alg:path-tracking} to track $\nabla_{x,\lambda} \mathcal{L}$'s solution from $(x_0, \lambda_0, 0)$ to $(x_\tau,\lambda_\tau,\tau)$, recording $x(t)$'s value at each step $0=t_0<t_1<\dots < t_N=T$ with $t_i=t_{i-1}+T/N$.}
\RETURN the discretized curve $x(t)$.
\end{algorithmic}
\end{algorithm}

\section{Deformation results for the tensegrity structures}
\label{section:modellingbmnsgn}

With a system of constraints for our tensegrity and a robust optimization algorithm at hand, we are finally able to tackle the problem of modeling the mechanical behavior of the tensegrity structures related to the $\Pi^+$ and $\Sigma^+$ rod packings from \cite{Evans2011}, which both display dilatant behavior. To set up the experiments, we first factor out rigid motions and consider the first coordinate as a parameter that will induce the frameworks' extension. Denoting the parameter by $\tau$, this leaves us with the lattice
\[
\Lambda_\tau = \begin{pmatrix}
    \tau & 0 & 0\\
    c_{12}(\tau) & c_{22}(\tau) & 0\\
    c_{13}(\tau) & c_{23}(\tau) & c_{33}(\tau)
\end{pmatrix}.
\]

\begin{figure}[h!]
    \centering
    \includegraphics[width=0.8\linewidth]{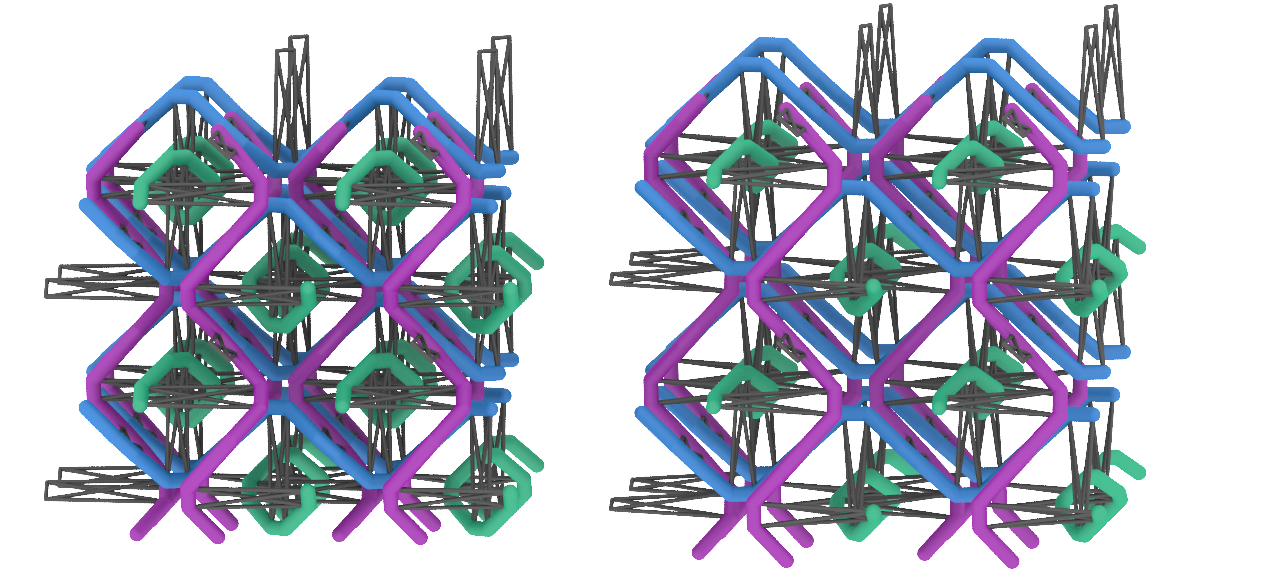}
    \vspace*{2mm}
    
    \includegraphics[width=0.8\linewidth]{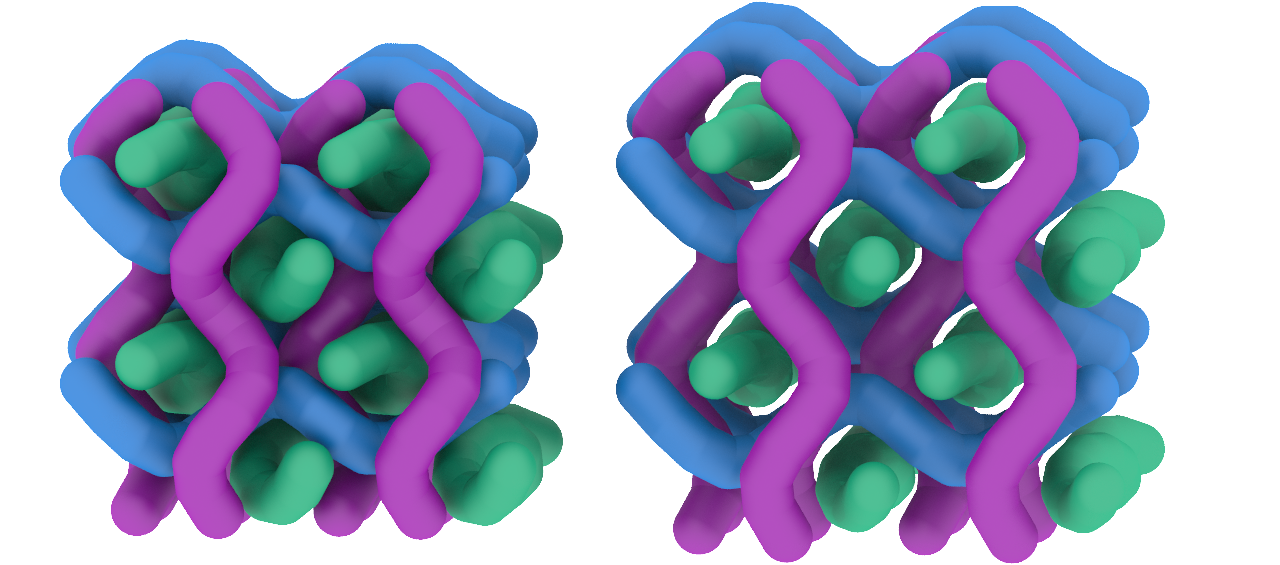}  
    \caption{\small Eight unit cells of the $\Pi^+$ tensegrity (t.) and its related cylinder packing (b.) at parameter values $\tau=1.0$ (l.) and $\tau=1.5$ (r.). The filaments' radii are displayed significantly smaller, without contacts, to better show how the structure deforms. The chosen projection is onto the $y$-$z$-plane, so the deformation's direction comes out of the page.}
    \label{fig:bmnanimation}
\end{figure}

In the case of the $\Pi^+$ tensegrity depicted in \Cref{fig:bmnanimation}, we record the framework's behavior for $\tau \in [0.93, 1.52]$, which amounts to an extension of roughly $63\%$. 
Beyond that, the structure in unstable. Having access to the lattice $\Lambda_\tau$, we can discretize the interval $\mathcal{I}_{\Pi^+} = [0.93,1.52]$ with step size $10 ^{-3}$ and perform the path-tracking \Cref{alg:path-tracking} to learn the lattice generators' value at each step. For the Poisson's ratio, just the orthogonal directions are relevant, so only $\Lambda_\tau$'s diagonal entries are involved in the calculation.

\begin{figure}[h!]
    \centering
    \includegraphics[width=0.93\linewidth]{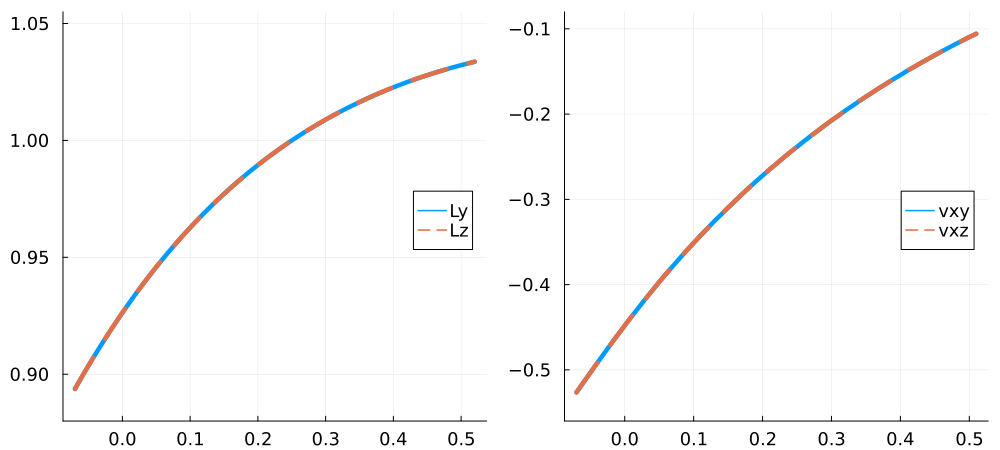}
    \caption{\small The lattice extensions $L_y$ and $L_z$ (l.) and the Poisson's ratios $\nu_{xy}$ and $\nu_{xz}$ (r.) corresponding to the homotopy of the $\Pi^+$ tensegrity, taking the model framework from $\tau=0.93$ to $\tau=1.52$. For each step with length $10^{-3}$, a data point is recorded. It can be observed that the lattice continuously grows from an extension of $0.89$ to $1.04$ in both $y$- and $z$-direction, amounting to a relative extension of $17\%$, compared to $63\%$ in $x$-direction. The Poisson's ratio ranges from $-0.52$ to $-0.11$, so the deformation path is auxetic.}
    \label{fig:bmnpoissonsratio}
\end{figure}

In \Cref{fig:bmnpoissonsratio}, the (orthogonal) lattice extensions $L_y$ and $L_z$ in $y$- and $z$-direction respectively are depicted, along with the corresponding Poisson's ratios $\nu_{xy}$ and $\nu_{xz}$. It allows us to conclude that the $\Pi^+$ tensegrity is auxetic in terms of the definition from materials science (cf. \Cref{section:previousmodel}), since the Poisson's ratio remains negative throughout the entire deformation path. The curves corresponding to the $y$- and $z$-direction are nearly identical, suggesting that the framework stays symmetric in these directions. \Cref{fig:bmnanimation} then depicts the tensegrity structures corresponding to two selected parameters from the chosen interval. We can again use the pictures to qualitatively deduce that the $\Pi^+$ packing extends in both $y$- and $z$-direction when stretched in $x$-direction. 

\begin{figure}[h!]
    \centering
    \includegraphics[width=0.8\linewidth]{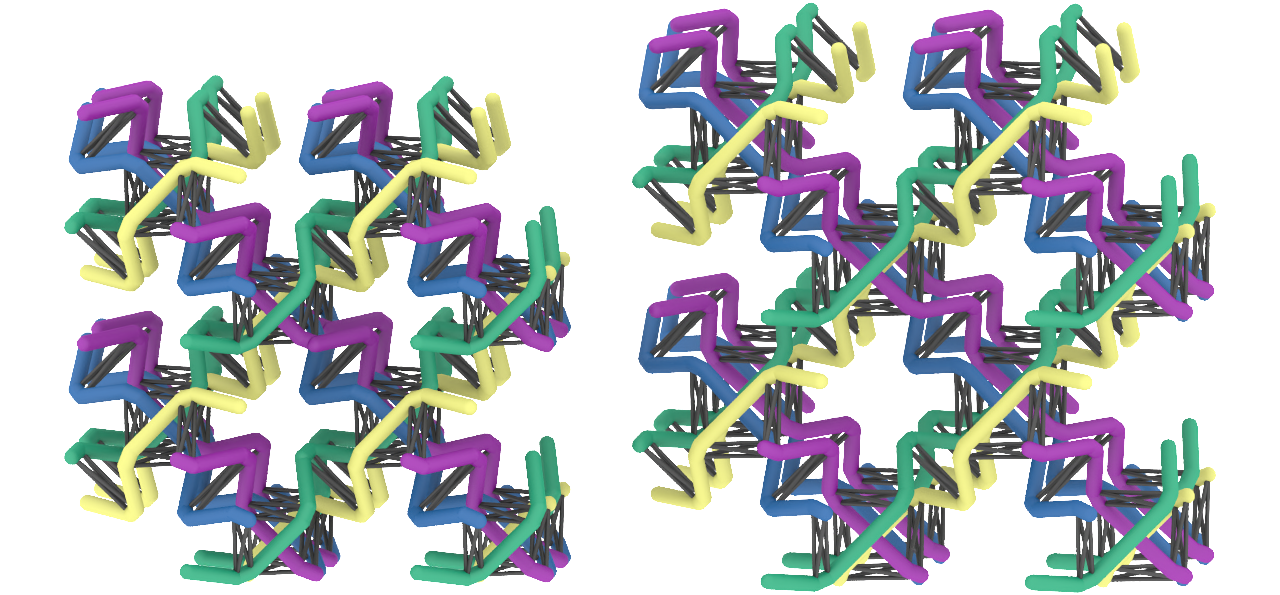}
    \vspace*{3mm}
    
    \includegraphics[width=0.8\linewidth]{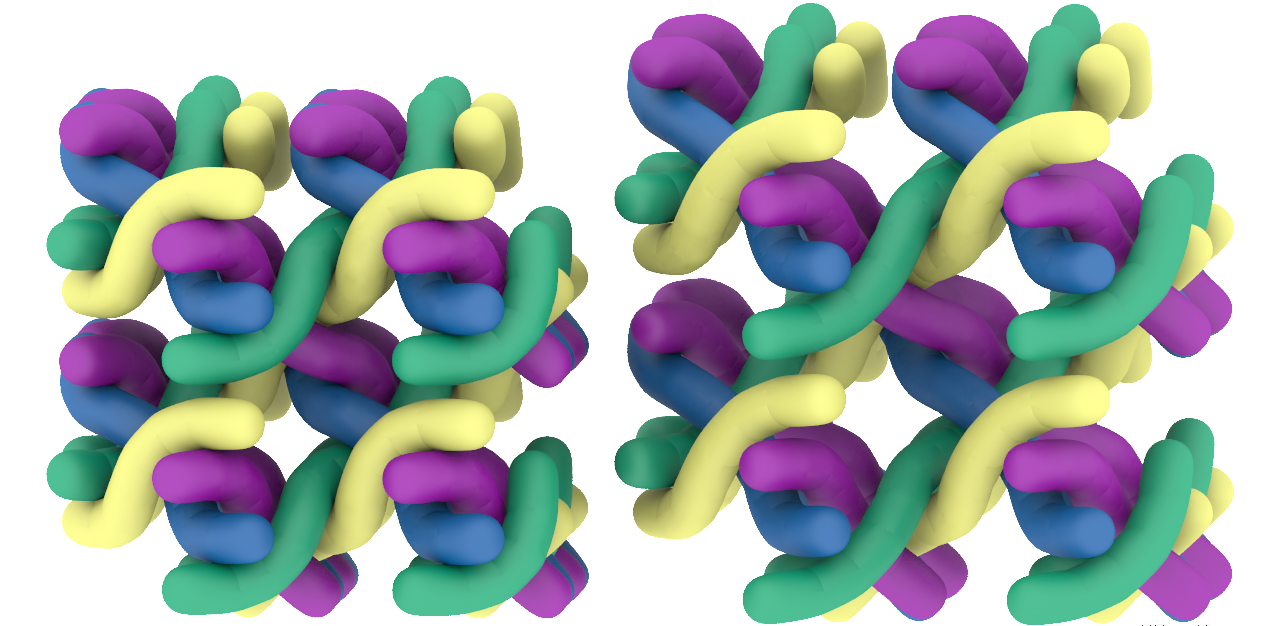}
    \caption{\small Eight unit cells of the $\Sigma^+$ tensegrity (t.) and its related cylinder packing (b.) at parameter values $\tau=1.0$ (l.) and $\tau=1.5$ (r.). The filaments' radii are displayed significantly smaller, without contacts, to better show how the structure deforms. The chosen projection is onto the $y$-$z$-plane, so the deformation is performed out of the page.}
    \label{fig:sgnanimation}
\end{figure}

\begin{figure}[h!]
    \centering
    \includegraphics[width=0.93\linewidth]{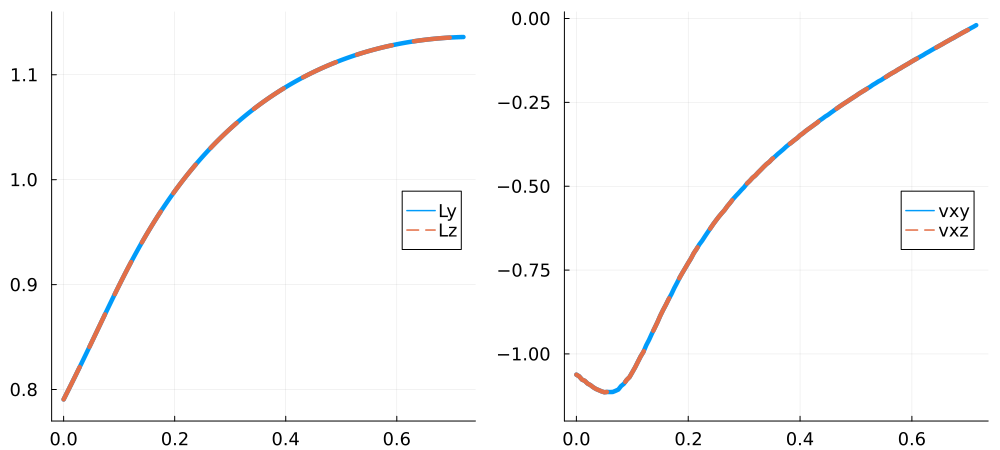}
    \caption{\small The lattice extensions $L_y$ and $L_z$ (l.) and the Poisson's ratios $\nu_{xy}$ and $\nu_{xz}$ (r.) corresponding to the homotopy of the $\Sigma^+$ rod packing, taking the model framework from $\tau=0.81$ to $\tau=1.52$. For each step with length $10^{-3}$, a data point is recorded. It can be observed that the lattice continuously grows from an extension of $0.79$ to $1.14$ in both $y$- and $z$-direction, amounting to a relative extension of $44\%$, compared to $88\%$ in $x$-direction. The Poisson's ratio grows from $-1.1$ over $-1.3$ until $-0.02$, so the deformation path is auxetic in the engineer's sense.}
    \label{fig:sgnpoissonsratio}
\end{figure}

We also analyze the $\Sigma^+$ tensegrity -- with two selected configurations depicted in \Cref{fig:sgnanimation}--, which is related to the structure found in skin cells. Similar to $\Pi^+$, we stretch the framework in the $x$-direction, recording the framework's behavior for $\tau \in [0.81, 1.52] = \mathcal{I}_{\Sigma^+}$. Again, the structure is unstable outside of that region. Discretizing the interval $\mathcal{I}_{\Sigma^+}$ with step size $10 ^{-3}$ and performing the path-tracking \Cref{alg:path-tracking} allows us to calculate the Poisson's ratio. In \Cref{fig:sgnpoissonsratio}, the (orthogonal) lattice extensions $L_y$ and $L_z$ along with the corresponding Poisson's ratios $\nu_{xy}$ and $\nu_{xz}$ are displayed.

The Poisson's ratio's curve is smooth and both curves are almost identical, ensuring the rod packing's symmetry in $y$- and $z$-direction. Still, the Poisson's ratio is not $-1$ everywhere, so the packing cannot extend perfectly symmetrically. Furthermore, the Poisson's ratio is negative throughout the entire deformation path, implying that auxeticity in terms of the definition from materials science is guaranteed.

What remains is to investigate, whether the tensegrity structures are also auxetic by the geometric interpretation given in \Cref{def:auxeticpath}. As was already proven in \Cref{prop:auxeticityandpoissonratio}, the geometric definition of auxeticity is stronger than just considering the Poisson's ratio. While the lattice $\Lambda_\tau$ is not exactly diagonal, throughout the deformation paths the off-diagonal entries are orders of magnitude smaller than the lattice's diagonal entries for the tetrahedral models associated to both $\Pi^+$ and $\Sigma^+$, giving us hope that the deformation path is geometrically auxetic after all.

Indeed, we can apply \Cref{prop:inequalitiesforauxeticity} to suggest that both tensegrity structures are auxetic.
To do so, let us start with the $\Pi^+$ tensegrity with corresponding deformation path parametrized by the interval $\mathcal{I}_{\Pi^+}$. First, we choose the step size $3 \cdot 10^{-3}$. 
Taking a look at the results our deformation path \Cref{alg:rodpacking} produced, we find that $c_{ii}(\tau_1)/c_{ii}(\tau_2)$ is bounded above by $0.99975$. 
With this, we can show that all three inequalities from \Cref{prop:inequalitiesforauxeticity} are satisfied for $\alpha= 10^{-4}$. We can then easily check that the off-diagonal entries' norms are bounded above by $10^{-5}$ throughout the deformation path, implying that the operator defined in Equation (\ref{eqn:linearoperatorauxetic}) is a contraction. 
\begin{figure}[h!]
    \centering
    \includegraphics[width=0.95\linewidth]{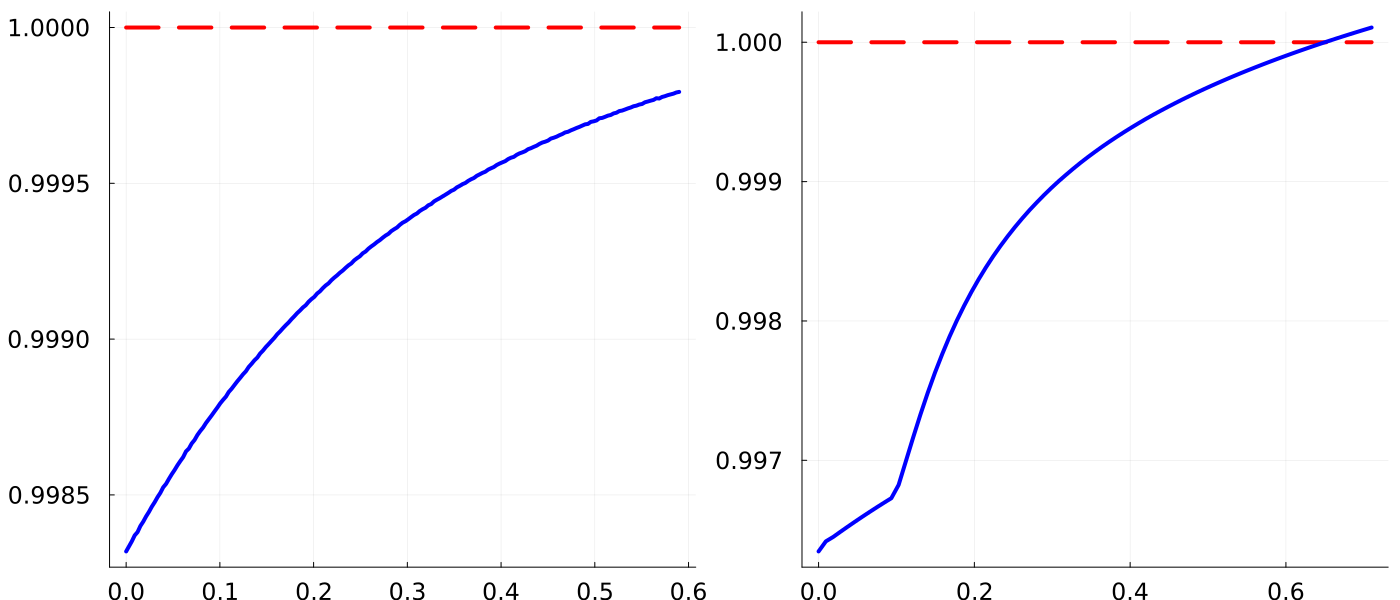}
    \caption{\small These diagrams depict the operator norms of the linear operator $T_{\tau_2\tau_1}$ taking the lattice $\pi_\tau$ corresponding to $\Pi^+$ tensegrity (l.) and the $\Sigma^+$ tensegrity (r.) from time $\tau_2$ to $\tau_1$ (see \Cref{section:geometricauxetics}) for fixed step length $\tau_2-\tau_1=3\cdot 10^{-3}$. Throughout most of the paths, the operator norm remains smaller than $1$, suggesting that the deformation paths are auxetic by \Cref{def:auxeticpath}, as expected. $\Sigma^+$ already becomes non-auxetic geometrically before the Poisson's ratio becomes positive.}
    \label{fig:operatornorm}
\end{figure}

By continuity, we can expect that this property persists for arbitrary discretizations and that the deformation path is actually auxetic as per \Cref{def:auxeticpath}. This result is documented by \Cref{fig:operatornorm}, which shows that the linear operator's norm is indeed bounded above by $1$.

 For the $\Sigma^+$ tensegrity, the deformation path is parametrized by the interval $\mathcal{I}_{\Sigma^+}$ with step size $3 \cdot 10^{-3}$. The numerical path-tracking \Cref{alg:path-tracking} returns $c_{ii}(\tau_1)/c_{ii}(\tau_2)\leq 0.999962$ and that the off-diagonal entries' norms are bounded above by $2\cdot 10^{-6}$. As \Cref{prop:inequalitiesforauxeticity} suggests, $\alpha = 1.5\cdot 10^{-5}$ suffices, implying that the operator defined in Equation (\ref{eqn:linearoperatorauxetic}) is a contraction. Again, continuity suggests that arbitrary discretizations will have this property, so the deformation path is likely auxetic in the geometric sense, too. \Cref{fig:operatornorm} verifies this result. Close to $\tau=1.55$, the Poisson's ratio is close to zero, which is consistent with the operator norms approaching 1.

\section{Conclusion}
\label{section:outlook}

We have explored the deformation of two examples of 3-periodic tensegrity structures, showing that they are auxetic. We analyzed the frameworks from a numerical optimization perspective, using the Euclidean distance retraction, as well as through existing techniques on geometric auxeticity. The result of this was a robust measurement of the Poisson's ratio of the tensegrity structures, which are clearly auxetic with well-behaved numerics. These results have a meaningful impact in theoretical materials science, where these structures could be potential targets for auxetic material design. They also demonstrate the usefulness of these numerical techniques on complicated 3-periodic framework materials, where a more systematic study of material properties is necessary.

The tensegrity structure was designed to mimic the mechanics of filament packings that showed a dilatent property related to auxeticity. The parallel between the filament packing mechanics and the tensegrity is still a work in progress, but we feel that this study is a good approximation. One complication is that the contacts of the filament packing are approximated by something designed for a perpendicular contact, but in a material this is not always the case. In general, we should get a good approximation in a material where the filaments contact each other roughly perpendicularly (and stay so throughout the deformation), and more problematic results when this is not the case. This is indeed what we see in further examples. Refining this model, and more accurately modeling the mechanics of filamentous materials, is the subject of future research. 

In summary, we have demonstrated the use of these mathematical techniques in obtaining meaningful results in materials science, where the high complexity of the structure is challenging for existing approaches. 
For the analysis of complicated material microstructures, we believe these techniques, which require research from both the materials science and mathematical perspective, are fruitful.

 {\footnotesize\linespread{0.8}

 \section*{Acknowledgments}
We want to thank Alexander Heaton for his collaboration on the the development of the package \texttt{HomotopyOpt.jl} and continued discussions about nonlinear optimization problems and robust path-tracking. We also thank Paul Breiding and Sascha Timme for making the package \texttt{HomotopyContinuation.jl}, that played a major role in the algorithm's development, available. 
 
\bibliographystyle{siamplain}
\bibliography{references}
}

\end{document}